\documentclass{amsart}
\usepackage{graphicx}
\usepackage{amssymb}
\usepackage{amsthm}
\usepackage{mathrsfs,txfonts}
\usepackage[all]{xy}
\usepackage{enumerate}
\DeclareMathAlphabet{\mathpzc}{OT1}{pzc}{m}{it}

\newcommand{\til}[1]{\widetilde{#1}}

\def\spec{\operatorname{Spec}}
\newcommand{\Oh}{\mathcal{O}}
\newcommand{\pf}{\noindent{\bf Proof:\ \ }}
\newcommand{\cqd}{{\hfill $\rule{2mm}{2mm}$}\vspace{3mm}}
\newcommand{\tmu}{\til{\mu}}

\def\mult{\operatorname{mult}}
\def\min{\operatorname{min}}

\def\deg{\operatorname{deg}}
\def\car{\operatorname{char}}

\def\lim{\operatorname{lim}}

\def\C{\mathbb{C}}
\def\O{{\mathcal{O}}}

\def\A{{\mathbb{A}}}

\def\N{{\mathcal{N}}}

\def\M{{\mathfrak{m}}}
\def\R{{\mathcal{R}}}

\def\mult{{\rm mult}}

\providecommand{\deg}{\mathop{\rm deg}\nolimits}

\providecommand{\ord}{\mathop{\rm ord}\nolimits}

\newtheorem{theorem}{Theorem}[section]
\newtheorem{lemma}[theorem]{Lemma}     
\newtheorem{corollary}[theorem]{Corollary}
\newtheorem{definition}[theorem]{Definition}
\newtheorem{proposition}[theorem]{Proposition}
\newtheorem{example}[theorem]{Example}
\newtheorem{remark}[theorem]{Remark}

\begin{document}
\title{The Milnor number of a hypersurface singularity in arbitrary characteristic}
 \author{A. Hefez, \ J.H.O. Rodrigues \ and \ R. Salom\~ao } 
\maketitle \vspace{-.7cm}

\begin{center} (Universidade Federal Fluminense - Niter\'oi) \end{center}
\vspace{1cm}

The Milnor number of an isolated hypersurface singularity, defined as the codimension $\mu(f)$ of the ideal generated by the partial derivatives of a power series $f$ whose zeros represent locally the hypersurface, is an important topological invariant of the singularity over the complex numbers, but its meaning changes dramatically when the base field is arbitrary. It turns out that if the ground field is of positive characteristic, this number is not even invariant under contact equivalence of the local equation $f$. In this paper we study the variation of the Milnor number in the contact class of $f$, giving necessary and sufficient conditions for its invariance. We also relate, for an isolated
sin\-gu\-la\-ri\-ty, the finiteness of $\mu(f)$ to the smoothness of the generic fiber $f=s$. Finally, we prove that the Milnor number coincides with the conductor of a plane branch when the characteristic does not divide any of the minimal generators of its semigroup of values, showing in particular that this is a sufficient (but not necessary) condition for the invariance of the Milnor number in the whole equisingularity class of $f$. \medskip

\medskip


\section{Introduction}
\label{I}

Let $\R=k[[X_1,\ldots,X_n]]$ be the ring of formal power series in $n$
indeterminates over an algebraically closed field $k$. We denote by $\M$ its maximal ideal and
by $\R^*$ its group of units. When $n=2$, we will write  $\R=k[[X,Y]]$.

Let $f\in \M \setminus \{0\}$. We define the hypersurface determined by $f$ as its class 
 under the associate equivalence relation: $(f)=\{ uf; \ u\in \R^*\}$.

Since we are interested in the geometry of the hypersurface $(f)$, we are led to consider the
quotient algebra $\O_f=\R/\langle f \rangle$, called the {\em local algebra} of
$(f)$.

We will consider two other equivalence relations in $\R$. 

$f$ and $g$ are said {\em right equivalent}, or {\em $r$-equivalent}, writing $f\sim_{r} g$, if there exists an automorphism $\varphi$ of $\R$ such that $g=\varphi(f)$. On the other hand, they are said {\em contact equivalent}, or {\em $c$-equivalent}, writing
$f \sim_{c} g$, if there exist an automorphism $\varphi$ and a unit $u$ of $\R$ such that
$ug=\varphi(f)$.

Notice that one has
\[
f\sim_{c} g \ \Longleftrightarrow \ \O_f \simeq  \O_g  \ \ \hbox{as $k$-algebras}.
\]

The {\em Jacobian ideal} of $f$ is the ideal generated by all the partial derivatives of
$f$
\[
J(f)=\langle f_{X_1}, \ldots, f_{X_n} \rangle.
\]
The {\em Milnor Algebra} of $f$ is the algebra $\R/J(f)$ and its
dimension $\mu(f)$ as a $k$-vector space is the {\em Milnor number} of $f$. If $f \sim_r g$,
we have that $\mu(f)=\mu(g)$.

The {\em Tjurina ideal} of $f$ is the ideal
\[
T(f)= \langle f \rangle + J(f)=\langle f, f_{X_1}, \ldots, f_{X_n} \rangle.
\]
The {\em Tjurina Algebra} of $f$ is the algebra $\R/T(f)$ and its
dimension $\tau(f)$ as a $k$-vector space is the {\em Tjurina number} of $f$.

Notice that when $\R=\C\{X_1,\ldots,X_n\}$ is the convergent power series ring, Milnor proves by
topological methods that $\mu(uf)=\mu(f)$ for any $u\in \R^*$. This result is usually
extended over characteristic zero fields by Lefschetz' principle. In arbitrary
characteristic this does not hold as one can see from the example below.

\begin{example} Let $\car{k}=p$ and $f=Y^p+X^{p+1}\in k[[X,Y]]$. Then $(f)$ is an irreducible curve such that $\mu(f)=\infty$, $\tau(f)=p^2$ and $$\mu((1+Y)f)=p^2 \neq \mu(f).$$
\end{example}

This is a relevant issue in our investigation since some important problems are
connected to it. For instance, we will characterize, in Section 3, the $f$ for which
$\mu(f)=\mu(uf)$ for all $u\in \R^*$ and will study in general the variation of $\mu(uf)$ when
$f$ is fixed and $u$ varies in $\R^*$.\medskip

From the inclusion $J(f) \subset T(f)$, it is clear that
\[\tau(f) \leqslant \mu(f).
\]

So, one always has
\[
\mu(f) < \infty \ \Rightarrow \tau(f) < \infty.
\]

In characteristic zero, one also has the converse of the above implication. This will be
proved algebraically in Section 2. In positive characteristic, the converse may fail,
as one can see from Example 1.1 above. \medskip

Motivated by the above discussion and by the fact that the ideal of a singularity on a hypersurface is the Tjurina ideal, the natural definition for isolated singularity is
the following:
\begin{definition} A hypersurface $(f)$ has an isolated singularity
at the origin if $f\in \M^2$ and $\tau(f)<\infty$.
\end{definition}

Notice that this is a well posed definition, since $\tau(f)=\tau(g)$ when $f$ and $g$
are contact equivalent. So, in characteristic zero, to say that $(f)$ has an isolated
singularity is equivalent to say that $\mu(f)<\infty$, but not in arbitrary
characteristic.

There is an easy criterion in arbitrary characteristic (cf. \cite{B}, Proposition 1.2.11)
for a plane curve $(f)$ to have an isolated singularity: $(f)$
has an isolated singularity if and only if $f$ is reduced. In contrast, the fact that
$f$ is reduced is not sufficient to guarantee that $\mu(f)<\infty$ as shows Example 1.1.
Also, the vanishing of one of the partial derivatives of $f$ implies $\mu(f)=\infty$,
but this is not a necessary condition, even in the case of plane curves, as the
following example shows.

\begin{example} Let $\car{k}=3$ and $f=X^2Y+Y^2X\in k[[X,Y]]$. We have that $f=XY(X+Y)$ is reduced, but $f_X$ and $f_Y$ have the common factor $Y-X$, implying that $\mu(f)=\infty$.
\end{example}

In the following section we will give a criterion for the finiteness of $\mu(f)$ which
will shed light on why in characteristic zero one has $\tau(f)< \infty$ implies $\mu(f)<\infty$.

Recall that in the complex case the Milnor number of $f$ was introduced in \cite{Mi} as the rank of the middle cohomology group of the fiber of the local smoothing $f=s$. In this
setting $\mu(f)$ is referred to as the number of vanishing cycles associated to $f$.
When we switch to a field of positive characteristic $p$, the fibration $f=s$ may not be a
local smoothing anymore, that is, it may be a counter example to Bertini's theorem on
the variation of singular points in linear systems, true in characteristic zero. In Section 4, we characterize this phenomenon, that may only occur in positive characteristic, in terms of the infiniteness of the Milnor number.

Finally, in the last two sections, we study plane branches singularities over arbitrary algebraically closed fields. In characteristic zero, the Milnor number $\mu(f)$ coincides with the conductor $c(f)$
of the semigroup of values of a branch $(f)$. In arbitrary characteristic, Deligne proves in \cite{De} (see also \cite{MH-W}) the inequality $\mu(f)\geqslant c(f)$, where the difference $\mu(f)-c(f)$ measures the
existence of wild vanishing cycles. We prove that Milnor's number and the conductor of a branch $(f)$ coincide when the characteristic does not divide any of the minimal generators of the semigroup of values of $f$. Our proof was inspired by a result of P. Javorski in the work \cite{Ja2}, which we simplified and extended to arbitrary characteristic, under the appropriate assumptions. We would like to point out that in the process of writing the final version of this paper, E. Garc\'ia Barroso and A. Ploski posted the preprint \cite{GB-P}, showing by other methods our last result (with the converse), but in the particular case when $p$ is greater than the multiplicity of $f$, and also observed that their proof fails when $p$ is less or equal than the multiplicity of $f$. We should also mention that H.D. Nguyen in \cite{Ng} has shown, in the irreducible case, the weaker result, namely, that if $p> c(f)+\mult(f)-1$, then $\mu(f)=c(f)$. Notice that once fixed the ground field $k$ of positive characteristic $p$, both results in \cite{GB-P} and \cite{Ng} cover only finitely many values of the multiplicity $\mult(f)$, while our result is in full generality.

This work is part of the PhD Thesis of the second author, under the supervision of the other two authors.



\section{The finiteness of $\mu(f)$}

To discuss the finiteness of $\mu(f)$ we must impose that $\tau(f)<\infty$, or
equivalently that $(f)$ has an isolated singularity, because otherwise $\mu(f)=\infty$.

The following result gives a criterion for the finiteness of $\mu(f)$.

\begin{proposition} \label{finiteness} Let $f\in \M$ and suppose that $\tau(f)<\infty$. Then
\[
\mu(f)<\infty \ \Leftrightarrow \ f\in \sqrt{J(f)}.
\]
\end{proposition}
\pf Notice that when $f\in \M \setminus \M^2$ both conditions hold trivially. So, we are
concerned with the case $f\in \M^2$.

Suppose that $\mu(f)<\infty$, then $J(f)$ is $\M$-primary, hence
$f\in \M=\sqrt{J(f)}$.

Conversely, suppose that $f\in \sqrt{J(f)}$. Since $\tau(f)<\infty$, we have that
$\sqrt{T(f)} =\M$. Now,
\[
T(f)=\langle f \rangle + J(f) \subset \sqrt{J(f)} \subset \M,
\]
since $f\in \M^2$. Taking radicals we get
\[
\M =\sqrt{T(f)} \subset \sqrt{J(f)} \subset \M.
\]
So, $\sqrt{J(f)} = \M$, which in turn implies that $\mu(f)<\infty$.
\cqd
 
\begin{remark} \ The reason why in characteristic zero
$\tau(f)<\infty \Rightarrow \mu(f)<\infty$ is that in this case one has
\[
f\in \overline{ \langle X_1f_{X_1},\ldots, X_nf_{X_n}\rangle },
\]
where the notation $\overline{I}$ means the integral closure of the ideal $I$
(cf. \cite{H-S}, Theorem 7.1.5). But, the above inclusion implies that
\[
f\in \overline{\M\, J(f)} \subset \overline{J(f)} \subset \sqrt{J(f)},
\]
and because of Proposition 2.1, $\tau(f)<\infty$ implies that $\mu(f)<\infty$.

Notice that the condition $f\in \sqrt{J(f)}$, which appears in Proposition 2.1 is
weaker than the condition $f\in \overline{J(f)}$ that holds in characteristic
zero. On the other hand, the condition $f\in \overline{J(f)}$, in arbitrary
characteristic, implies the Brian\c con-Skoda inclusion:  $f^n\in J(f)$, where
$n=\dim \R$ (cf. \cite{H-S}, Theorem 13.3.3).
\end{remark}

Recall that an ideal $J$ is called a reduction of an ideal $T$ if $J\subset T$ and
there exists $n\in \mathbb N$ such that $T^{n+1}=JT^n$. We denote by $e_0(I)$ the
Hilbert-Samuel multiplicity of an $\M$-primary ideal $I$ and put $e_0(I)=\infty$ if
$\sqrt{I}\subsetneq \M$.

The next well known proposition describes, in general, the effect of the condition
$f\in \overline{J(f)}$.

\begin{proposition}[\cite{N-R} and \cite{R}] \label{NR}
Let $f\in \M$ be such that $\tau(f)<\infty$. The following conditions are
equivalent:\smallskip
\begin{enumerate}[{\rm (i)}]
\item $f\in \overline{J(f)}$;
\item $J(f)$ is a reduction of $T(f)$;
\item  $e_0(T(f))=e_0(J(f))=\mu(f)$.
\end{enumerate}
\end{proposition}
\pf (i) $\Leftrightarrow$ (ii): cf. [H-S, Corollary 1.2.2].

\noindent (ii) $\Leftrightarrow$ (iii): This follows from \cite{N-R} and from
\cite{R}, Theorem 3.2, since $\R$ is a level ring (analytically unramified).
\cqd

\begin{corollary}Let $k$ be a field of characteristic zero and
$f\in\mathfrak{m}$. Then $\mu(f)$ is invariant under contact equivalence.
\end{corollary}

\textit{Proof:} If $\tau(f)=\infty$ then $\mu(g)=\infty$, for every $g$ in the same
contact equivalence class and we are done. So we are restricted to the case
$\tau(f)<\infty.$ Since changing coordinates obviously does not change $\mu$, we only
need to show that $\mu(f)=\mu(uf)$ for every unit $u\in \R$.
However it is easy to see that $T(f)=T(uf)$, for every such $u$. In characteristic zero,
both $J(f)$ and $J(uf)$ are reductions of $T(f)$, according to the previous proposition
and remark. On the other hand, $\R$ is a regular (hence Cohen-Macaulay) local ring. So
$\mu(f)=e_0(J(f))$ (cf. \cite{Ma}, Theorem 17.11). Therefore,

\[\mu(f)=e_0(J(f))=e_0(T(f))=e_0(T(uf))=e_0(J(uf))=\mu(uf).\]\cqd

\begin{remark} From the inclusion $J(f)\subset T(f)$ one has that
$e_0(T(f)) \leqslant e_0(J(f))$ (cf. \cite{Ma}, Formula 14.4), then
\[
e_0(T(f))\leqslant e_0(J(f))=\mu(f).
\]
\end{remark}

The inequality in the above remark may be strict, as shows the following example.

\begin{example}
Let $\car{k}=p$ and $f=X^p+X^{p+2}+Y^{p+2}\in k[[X,Y]]$. 
As $J(f)=\langle X^{p+1},Y^{p+1} \rangle$, we have that $\mu(f)=(p+1)^2$ and
$\tau(f)=p(p+1)<\infty$. If $g=(1+X)f$, then
an easy calculation with intersection indices shows that $\mu(g)=I(g_X,g_Y)=p(p+1)$, so
\[
e_0(T(f))=e_0(T(g))\leqslant \mu(g) =p(p+1) < (p+1)^2=\mu(f).
\]
\end{example}

It follows from the preceding discussion that the importance of the Jacobian ideal of a hypersurface singularity in characteristic zero is due to the fact that it is a minimal reduction of the Tjurina ideal, which is the ideal that carries all the information about the singular point. In this situation, one has that $e_0(T(f))=\mu(f)$, and this is why Milnor's number is full of meanings in characteristic zero. 

This leads us to consider the {\it Milnor number of a hypersurface}
$(f)$ as
\[ \tmu(f)=e_0(T(f)), \]
which is an invariant of the contact class of $f$. Notice that in characteristic zero one always has $\tmu(f)=\mu(f).$

\begin{remark}
Proposition \ref{NR} gives a numerical criterion for testing if
$f$ belongs to $\overline{J(f)}$. Example 2.6 shows that one may have $f\in \sqrt{J(f)}$ with 
$f \not\in \overline{J(f)}$, since in this case $\mu(f)>e_0(T(f))$.
\end{remark}




\section{Variation of $\mu(uf)$ and computation of $\tmu(f)$}
 
We have seen that in characteristic zero the multiplicity $e_0(T(f))$ may be computed as
the codimension of $J(f)$ in $\R$ because, in that situation, $J(f)$ is a minimal
reduction of $T(f)$. On the other hand, this is not always the case if the ground field has positive characteristic. Therefore, we are led to investigate whether $J(f)$ is a minimal reduction of $T(f)$ when $\tau(f) < \infty$. As a consequence of our discussion we will analyze 
the variation of $\mu(uf)$ when $u$ varies in $\R^*$ and obtain a method for computing 
$\tmu(f)$. More generally, we will search for minimal reductions of $T(f)$. 

Since $\R$ is a local ring with infinite residue field $k$ and $f$ is such that $\tau(f)<\infty$, it is well known that for a fixed set of generators, not necessarily minimal,   $f,f_{X_1},\ldots,f_{X_n}$ of $T(f)$, if we take {\em sufficiently general} linear combinations
\begin{equation}\label{generallc}
g_i=h_{0,i}\, f+\sum_{j=1}^nh_{j,i}\; f_{X_j}, \ \ h_{j,i} \in \R, \ \ i=1,\ldots,n,
\end{equation}
then $g_1,\ldots,g_n$ is a system of parameters in $\R$ and the ideal they generate is
a reduction of the $\M$-primary ideal $T(f)$, hence a minimal reduction
(cf. \cite{Ma}, Theorem 14.14).

To find the conditions on the $h_{j,i}$ to be sufficiently general, we will need the notion
of null-forms.

A null-form (cf. \cite{Ma}, proof of Theorem 14.14) for the ideal $T(f)$ is a homogeneous
polynomial $\varphi \in k[Y_0,\ldots,Y_n]$ of some degree $s$ such that there exists
$F \in \R[Y_0,\ldots,Y_n]$ homogeneous of degree $s$ for which $F \equiv \varphi \bmod \M$ 
and $F(f,f_{X_1},\ldots,f_{X_n})\in \M \,T(f)^s$. This notion is independent of the choice of  $F$. We denote by $\N_{T(f)}$ the homogeneous ideal in $k[Y_0,\ldots,Y_n]$ generated by all  null-forms of $T(f)$.

\begin{remark} The ideal $\N_{T(f)}$ depends upon the generators
$f,f_{X_1}, \ldots, f_{X_n}$ of $T(f)$. As $k$-algebras one has
\[\dfrac{k[Y_0,\ldots,Y_n]}{\N_{T(f)}}\simeq\bigoplus_{s\geqslant 0}\dfrac{T(f)^s}{\mathfrak{m}
T(f)^s}.\]
The $k$-algebra on the right hand side is called the \textit{fiber cone}
of $T(f)$ and it is the graded ring corresponding to the special fiber of the blow-up of
$\spec \R$ at $T(f)$. We also have
\[
\dim_{krull} \dfrac{k[Y_0,\ldots,Y_n]}{\N_{T(f)}}= n
\]
(cf. \cite{Ma}, proof of Theorem 14.14), which implies that the projective zero set
$Z(\N_{T(f)})$ in ${\mathbb P}^{n}_k$ is of dimension $n-1$. In particular,
$\N_{T(f)}\neq \langle 0 \rangle$. \
\end{remark}

\begin{example} Recall Example 1.1, where $f=Y^p+X^{p+1}\in k[[X,Y]]$ and 
$\car{k}=p$. Since the polynomial $Y_2\in k[Y_0,Y_1,Y_2]$ vanishes when evaluated at $(f,f_X,f_Y)$, we have that $Y_2\in \N_{T(f)}$. So, $Z(\N_{T(f)})\subset Z(Y_2)$, and since
$\dim(Z(\N_{T(f)}))=\dim(Z(Y_2))$ and $Z(Y_2)$ is irreducible, we have that
$Z(\N_{T(f)})=Z(Y_2)$. But this last equality, together with $Y_2\in \N_{T(f)}$, imply
that $\N_{T(f)}=\langle Y_2 \rangle$. 
\end{example}

\begin{example} Going back to Example 1.3, if $f=X^2Y+Y^2X\in k[[X,Y]]$, where
$\car{k}=3$, we have that $Y_0(Y_1+Y_2)\in \N_{T(f)}$, since
\[
f(f_X+f_Y)=-XY(Y+X)(Y-X)^2=(Y+X)f_Xf_Y\in \M\,T(f)^2.
\]
Now, because $Y_0, Y_1+Y_2 \not \in \N_{T(f)}$, it follows that $\N_{T(f)}$ is not a
prime ideal.
\end{example}

Given $f\in \M$, in order to have $e_0(T(f))=\mu(uf)$, for some unit $u$, we must find $u\in \R^*$ such that $uf\in \overline{J(uf)}$. We will show next that this is so for general units.\medskip

We will need the following result due to Northcott and Rees
(\cite{N-R}, or \cite{Ma}, proof of Theorem 14.14):

{\em The ideal $\langle g_1,\ldots,g_n\rangle$, where the $g_i$'s are as in (\ref{generallc}) is a
reduction of the ideal $T(f)$ if and only if the linear forms
\[
\ell_i=\sum_{j=0}^n a_{j,i}Y_j, \ \ i=1,\ldots,n,
\]
where $a_{j,i} \equiv h_{j,i} \bmod \M$, are such that the ideal
$\N_{T(f)}+\langle \ell_1, \ldots,\ell_n \rangle$ is
$\langle Y_0,\ldots,Y_n \rangle$-primary, that is,
$Z(\N_{T(f)}+\langle \ell_1, \ldots,\ell_n \rangle)=\{0\}\subset
\mathbb{A}^{n+1}_k$.}

\begin{theorem} Let $f\in \M$ and $u=\alpha_0+\alpha_1X_1+\cdots+\alpha_nX_n+\, hot.$,
with $\alpha_0\neq 0$ be a unit in $\R^*$. We have that $uf\in \overline{J(uf)}$ if
and only if there exists $G\in \N_{T(f)}$ such that
$G(\alpha_0,-\alpha_1,\ldots,-\alpha_n)\neq 0$. In particular, this holds for a generic
$(\alpha_0:\cdots : \alpha_n)\in \mathbb{P}^n_k$.
\end{theorem}
\pf If $g=uf$, then
\[
g_{X_i}=u_{X_i}f+uf_{X_i}, \ \ i=1,\ldots,n,
\]
with associated linear forms
\[
\ell_i=\alpha_iY_0+\alpha_0Y_i,  \  \ i=1,\ldots,n.
\]

We then have
{\small \[
\begin{array}{c}
Z(\N_{T(f)}+\langle \ell_1,\ldots, \ell_n \rangle) =\\ \\

Z\left(\left\langle\, \left( \frac{Y_0}{\alpha_0}\right)^{\deg(G)}
G(\alpha_0,-\alpha_1,\ldots,-\alpha_n), \, \alpha_1 Y_0+\alpha_0Y_1,\ldots,
\alpha_n Y_0+\alpha_0Y_n ; \ G\in \N_{T(f)}\setminus \{0\}\right\rangle\right).
\end{array}
\]}

Since $uf\in \overline{J(uf)}$ if and only if
$J(uf)=\langle g_{X_1},\ldots,g_{X_n}\rangle$ is a reduction of $T(uf)=T(f)$, then 
from the Northcott and Rees Theorem above mentioned, this happens if and only if
$Z(\N_{T(f)}+\langle \ell_1,\ldots,\ell_n\rangle)=\{0\}$.
This, in turn, happens if and only if for some $G\in \N_{T(f)}$, one has
$G(\alpha_0,-\alpha_1,\ldots,-\alpha_n)\neq 0$. 
\cqd

The above theorem shows that if $u$ is a general unit, in the sense that it has a
general linear term, then $J(uf)$ is a reduction of $T(uf)$, and so,
$uf\in \overline{J(uf)}$, which in turn implies that
\[
\mu(uf)=e_0(T(uf))=\tmu(f).
\]

This theorem allows us to give the following interpretation for $\tmu(f)$. 

\begin{corollary} Let $f\in \M$ be such that $\tau(f)<\infty$. Then $\tmu(f)=\min\{ \mu(uf), \ u\in \R^* \}$ and $\mu(f)=\tmu(f)$ if and only if $f\in \overline{J(f)}$.
\end{corollary}
\pf From Remark 2.5 we know that $\tmu(f)=\tmu(uf)\leqslant \mu(uf)$, for all $u\in \R^*$. According to the last theorem there is a unit $v$ such that $vf\in \overline{J(vf)}$, hence $\tmu(f)=\mu(vf)$, proving the first assertion. The second assertion follows immediately from Proposition \ref{NR}. \cqd

We also have the following result.

\begin{theorem}\label{mustavel} Let $f\in \M$ be such that $\tau(f)<\infty$. The following three statements are
equivalent.
\begin{enumerate}[{\rm (i)}]
\item $\mu(uf)=\tmu(f)$ for every unit $u\in \R$;
\item $Z(\N_{T(f)})=Z(Y_0)$;
\item $f^{\ell}\in\mathfrak{m}T(f)^{\ell}$, for some $\ell\geqslant 1$.
\end{enumerate}
\end{theorem}
\pf
\textbf{$(i\Rightarrow ii)$} If $\mu(uf)=\tmu(f)$ for every unit $u$,
then $Z(\N_{T(f)})\cap \{Y_0\neq 0\}=\emptyset$. Otherwise, if
$(1:\beta_1:\cdots:\beta_n)$ is in this set, then we would have $G(1,\beta_1,\ldots,\beta_n)=0$, for every $G\in \N_{T(f)}$, hence $u=1-\beta_1X_1-\cdots-\beta_nX_n$ would be such that
$\mu(uf)>\tmu(f)$, contradicting the hypothesis. 

Therefore, $Z(\N_{T(f)})\subseteq Z(Y_0)$, which implies the equality by comparing dimensions and by the irreducibility of $Z(Y_0)$.

\noindent \textbf{$(ii \Rightarrow iii)$} If $Z(\N_{T(f)})=Z(Y_0)$ then
$\sqrt{\N_{T(f)}}=\sqrt{\langle Y_0\rangle}=\langle Y_0\rangle$. Hence, there exists a positive integer $\ell$
such that $G=Y_0^{\ell}\in N_{T(f)}$. In other words,
$f^{\ell}\in \mathfrak{m}T(f)^{\ell}$.

\noindent \textbf{$(iii\Rightarrow i)$} If for some $\ell$, one has $f^{\ell}\in \mathfrak{m}T(f)^{\ell}$,
then $G=Y_0^{\ell}\in N_{T(f)}\setminus \{0\}$. Let
$u=\alpha_0+\alpha_1X_1+\cdots+\alpha_nX_n+hot$ be a unit, then
$G(\alpha_0,-\alpha_1,\ldots,-\alpha_n)=\alpha_0^{\ell}\neq 0$. It follows that
$uf\in\overline{J(uf)}$ and, therefore, $\mu(uf)=\tmu(f)$.\cqd

\begin{corollary} Suppose $p=\car{k}=0$ and let $f \in \mathfrak{m}\setminus \{0\}$ with  $\tau(f)<\infty$. Then there exists $\ell\geqslant 1$ such that $f^{\ell}\in\mathfrak{m}T(f)^{\ell}.$
\end{corollary}

\pf For $p=0$ we know that $\tmu(f)=\mu(f)=\mu(uf)$ for every invertible $u$,
hence we may use the preceding theorem.\cqd

\begin{remark}[cf. \cite{Ga}] The preceding corollary can be derived from the
fact that if $p=0$, then $f\in\overline{\mathfrak{m}\, J(f)},$ as we have already
seen. Indeed, consider an equation of integral dependence of $f$ over
$\mathfrak{m}\, J(f)$:
\[f^{\ell}+a_1 f^{\ell-1}+\cdots+a_{\ell}=0,\]
with $a_i\in (\mathfrak{m}\, J(f))^i$. Hence, for each $i$,
\[f^{\ell-i}a_i\in\mathfrak{m}^i\,f^{\ell-i}\,J(f)^{i}\subset
\mathfrak{m}^i\,T(f)^{\ell}\subset\mathfrak{m}\,T(f)^{\ell}\]
and we conclude since $f^{\ell}=-\sum_i a_if^{\ell-i}.$ However, if $p>0$ one may produce an example of a power series $f$ which satisfies the equivalent conditions
of the Theorem \ref{mustavel}, but such that $f\not\in\overline{\mathfrak{m}\, J(f)}$.
\end{remark}

If $f$ is such that $\mu(uf)$ is independent of the unit $u$, that is, when $\mu(f)$
is invariant under contact equivalence, we will say that $f$ is {\em $\mu$-stable}.

The third condition in Theorem \ref{mustavel} may help to decide whether a given power
series is or not $\mu$-stable as we can see in the following examples. However, in order to
have this condition as a computational method, we need to bound $\ell$.

\begin{example} Let $f\in k[X_1,\ldots,X_n]$ be a quasi-homogeneous polynomial of degree $d$ and $\car{k}=p>0$. If $p \nmid d$, then $f$ is  $\mu$-stable.
Indeed, there are integers $d_1,\ldots,d_n$ such that
\[
df=d_1X_1f_{X_1}+\cdots+d_nX_nf_{X_n},
\]
which, since $p\nmid d$, implies that $f\in \M\, T(f)$, so by
Theorem \ref{mustavel}, $f$ is $\mu$-stable.
\end{example}

\begin{proposition} \label{multiplicity} Let $f=L^d+\,hot\;\in \R$, where $L$ is a linear form and $\car k=p$. If $f$ is $\mu$-stable, then $p\nmid d$.
\end{proposition}
\pf We will actually show that if $p|d$, then $f^{\ell}\not \in \M\, T(f)^{\ell}$,
for all $\ell \in \mathbb N$.

By a linear change of coordinates we may assume that $f=X_1^d+g$, where
$g\in \M^{d+1}$. It follows that $\mult(f_{X_1})\geqslant d-1$ and
$\mult(f_{X_i})\geqslant d$ for $i>1$. Now, we have
\[
T(f)^{\ell}=\langle f^{\alpha_0}f_{X_1}^{\alpha_1}\cdots f_{X_n}^{\alpha_n},
\ \alpha_0+\cdots+\alpha_n=\ell\rangle,
\]
and
\[
\mult(f^{\alpha_0}f_{X_1}^{\alpha_1}\cdots f_{X_n}^{\alpha_n}) \geqslant
\alpha_0d+\alpha_1\mult(f_{X_1})+(\alpha_2+\cdots+\alpha_n)d.
\]

So, if $p|d$, then $\mult(f_{X_1}) \geqslant d$, and therefore
\[
\mult(f^{\alpha_0}f_{X_1}^{\alpha_1}\cdots f_{X_n}^{\alpha_n})\geqslant
\sum_{i=0}^n\alpha_id=\ell d.
\]
This implies that $\mult(h) \geqslant \ell d+1$, for all $h\in \M\,T(f)^\ell $; and
since
$\mult(f^\ell )=\ell d$, it follows that $f^\ell \not \in \M\,T(f)^\ell$, for all
$\ell \in \mathbb N$.
\cqd

The above proposition has the following corollary:

\begin{corollary}
For $f\in k[[X,Y]]$ irreducible, where 
$\car k=p$, one has \[
f \ {\rm is} \ \mu-\hbox{stable} \ \Rightarrow \ p\nmid \mult(f).
\]
\end{corollary}

The converse of the above corollary  is not true, as one may see from the following example:

\begin{example}
Let $f=Y^3-X^{11}$ where $\car{k}=11$. Then $p\nmid \mult(f)$ but $f$ is not
$\mu$-stable, since $\mu(f)=\infty$ and $\mu((1+X)f)=22$.

Notice that whether $f$ is $\mu$-stable, or not, depends upon the
characteristic $p$ of the ground field. For example, when $p=5$ the same $f$ as above is quasi-homogeneous of degree $d=33$ which is not divisible by $p$. Hence it is $\mu$-stable.
\end{example}

\begin{remark} For all irreducible $f\in \M \setminus \M^2$, that is, for all irreducible smooth hypersurface germs, one has $\mu$-stability, since in such case $\tmu(f)=\mu(f)=0$.
\end{remark}

Finally, we give an example to show that the $\mu$-stability is not preserved by blowing-up.

\begin{example} Let $\car{k}=2$, and $f=Y^3+X^5\in k[[X,Y]]$. In this case,
we have that $f$ is $\mu$-stable, since $f\in \M\,T(f)$, but its strict transform
$f^{(1)}=v^3+u^2$ is not $\mu$-stable (cf. Proposition \ref{multiplicity}).
\end{example}




\section{Bertini's Theorem vs Milnor Number}

Let $f \in k[X_1,\ldots,X_n]$ with an isolated singularity at the origin of $\A^n$.
In this section we are going to study when the fibration $f\colon \mathbb{A}^n \to \mathbb{A}^1$
is a local smoothing of the singularity. Notice that in characteristic zero this is always the case, according to Bertini's theorem on the variation of singular points in linear systems.
However, it is well known that this is not true over fields of positive characteristic.

\begin{example}
 We have already seen that $f=X^p+Y^{p+1}$ has an isolated singularity at the origin
 when $\car{k}=p>0$. The fiber over each $s\in \mathbb{A}^1$ has $(s^{1/p},0)$ as a
 singularity.
 \end{example}
In positive characteristic we have the following characterization.

\begin{theorem}
Let $f$ be a polynomial admitting an isolated singularity at the origin.
 The fibration $f\colon \mathbb{A}^n\to \mathbb{A}^1$ is a local smoothing 
if and only if $\mu(f)<\infty$.
\end{theorem}
\pf If $\mu(f)=\infty$, then the codimension of the Jacobian ideal
$J(f)=\langle f_{X_1}, \ldots, f_{X_n}\rangle$ in $\mathcal{O}_{\mathbb{A}^n,0}$ is infinite.
This implies that the sequence $f_{X_1}, \ldots, f_{X_n}$ is not
$\mathcal{O}_{\mathbb{A}^n,0}$-regular. In this case $Z(f_{X_1}, \ldots, f_{X_n})$
contains a curve $C$ trough the origin in $\mathbb{A}^n$. We clearly
have that $C\cap Z(f-s)$ is a singular point of the fiber $f=s$.
Hence, it remains to show that $C$ dominates $\mathbb{A}^1$ under $f$. Otherwise,
$f(C)$ would be finite and there might exist $s_0$ such that
$Z(f_{X_1}, \ldots, f_{X_n},f-s_0)$ is infinite in some neighborhood at the origin
of $\mathbb{A}^n$. But this is a contradiction because, if $s_0\neq 0$, then $f-s_0$
does not vanish in some neighborhood of the origin and if $s_0=0$ it says that $f$ does not
have an isolated singularity at the origin.

Now, if $\mu(f)<\infty$ then the same argument used in Proposition \ref{finiteness}
shows that $f$ belongs to the ideal
$\sqrt{\langle J(f)\rangle}$ of
 $\Oh_{\A^n,0}$. Hence there exists a relation
\begin{equation}
\label{relation}
B\,f^N=A_1f_{X_1}+ \cdots + A_nf_{X_n}, \ \ \hbox{with} \ A_1,\ldots, A_n,B\in k[X_1,\ldots,X_n], \ B(0)\neq 0.
\end{equation}
Notice that each fiber $f^{-1}(s)$, with $s\neq 0$ in the open neighborhood $\A^n\setminus Z(B)$ of the origin, is smooth. Indeed, if $x \in \A^n  \setminus Z(B)$
 is a singular point of the fiber $f^{-1}(s)$, with $s\neq 0$, then
 $f(x)=s$ and $f_{X_i}(x)=0$ for each $i=1,\ldots, n$.
 On the other hand, since $B(x)\neq 0$ it follows from (\ref{relation})
 that $s=f(x)=0$, which is a contradiction. \cqd 

 
 

\section{Milnor number for plane branches with tame semigroups}

For the definitions and notation used in this section we refer to \cite{He}. Let $f \in \M \subset k[[X,Y]]$ be an irreducible power series, where $k$ is an algebraically closed field of characteristic $p\geqslant 0$. In this situation, we call the curve $(f)$ a {\em plane branch}. Let us denote by $S(f)=\langle v_0,\ldots,v_g\rangle$ the semigroup of values of the branch $(f)$, represented by its minimal set of generators. These semigroups have many special properties which we will use throughout this section and describe them briefly below.

Let us define $e_0=v_0$ and denote by $e_i=\gcd\{v_0,\ldots,v_{i}\}$ and by $n_i=e_{i-1}/e_i$,  $i=1,\ldots,g$. The semigroup $S(f)$ is {\em strongly increasing}, which means that $v_{i+1}>n_iv_i$, for $i=0,\ldots,g-1$, (cf. \cite{He}, (6.5)). This implies that the 
the sequence $v_0,\ldots,v_g$ is {\em nice}, which means that $n_iv_i\in \langle v_0,\ldots,v_{i-1}\rangle$, for $i=1,\ldots,g$, (cf. \cite{He}, Proposition 7.9). This, in turn, implies that the semigroup $S(f)$ has a {\em conductor}, denoted by $c(f)$, which is the integer characterized by the following property: $c(f)-1\not\in S(f)$ and $x\in S(f)$, for all $x\geqslant c(f)$, and it is given by the formula (cf. \cite{He}, (7.1))
\[
c(f)=\sum_{i=1}^g (n_i-1)v_i-v_0+1.
\]

The semigroup $S(f)$ is also symmetric (cf. \cite{He} Proposition 7.7), that is,
\[
\forall z\in \mathbb N, \ z\in S(f) \ \Longleftrightarrow \ c-1-z \not\in S(f).
\]

To deal with the positive characteristic situation, we introduce the following definition:

We call $S(f)$ a {\it tame semigroup} if $p$ does not
divide $v_i$ for all $i \in \{0,\ldots, g\}$.

Recall that two plane branches over the complex numbers are {\it equisingular} if their semigroups of values coincide. We will keep this terminology even in the case of positive characteristic.

The following example will show that $\tilde{\mu}(f)$ may be not constant in an
equisingularity class of plane branches. 

\begin{example} The curves given by $f=Y^3-X^{11}$ and $h=Y^3-X^{11}+X^8Y$ are
equisingular with semigroup of values $S=\langle 3,11\rangle$, but in characteristic $3$, one has $\tilde{\mu}(f)=\mu((1+Y)f)=30$, because $Y_2\in \N_{T(f)}$ and $\tilde{\mu}(h)=\mu((1+X)h)=24$, because $Y_1^3\in \N_{T(h)}$. Notice that in this case $S$ is not tame.
\end{example}

\begin{remark} For the above $h$ one has $\mu(h)=\infty$, $\tilde{\mu}(h)=24$ and $\tau(h)=22$. This shows that there is no isomorphism $\varphi$ of $k[[X,Y]]$ and no $H\in k[[X,Y]]$ such that
$\varphi(h)=H(X,Y^3)$, because, otherwise, we would get the contradiction
$$24=\tilde{\mu}(h)= \tilde{\mu}(\varphi(h))=\tilde{\mu}(H(X,Y^3))=\tau(H(X,Y^3))=\tau(h)=22.$$
\end{remark}

The following is an example which shows that the $\mu$-stability is not a character of an equisingularity class.

\begin{example} Let $S=\langle 4,6,25\rangle$ be a strongly increasing semigroup with conductor $c=28$. Consider the equisingularity class determined by $S$ over a field of characteristic $p=5$. If $f=(Y^2-X^3)^2-YX^{11}$, which belongs to this equisinsingularity class, we have that $\mu(f)=41$ and $\tilde{\mu}(f)=30$, hence $f$ is not $\mu$-stable. But the equisingular curve $h=(Y^2-X^3+X^2Y)^2-YX^{11}$ is $\mu$-stable, since $Y_0^3\in \N_{T(h)}$.
In this case one has $\mu(h)=\tilde{\mu}(h)=29$. Notice that here, again, $S$ is not tame. 
\end{example}

The aim of this section is to prove the following result:

\begin{theorem}[Main Theorem]\label{tame} If $f\in \M^2$ is a plane branch singularity with $S(f)$ tame, then  $\mu(f)=\tmu(f)=c(f)$. In particular, $f$ is $\mu$-stable.
\end{theorem}

The proof we give of this theorem is based on the following theorem which was stated without a proof over the complex numbers in \cite{Ja1}, but proved in the unpublished work \cite{Ja2}. Our proof, in arbitrary characteristic, is inspired by that work, which we suitably modified in order to make it work in the more general context we are considering.  

\begin{theorem} \label{theorem1}
Let $f\in \M^2$ be an irreducible Weierstrass polynomial such that $S(f)$ is tame. Then any family $\mathcal{F}$ of elements inside $k[[X]][Y]$ of degree in $Y$ less than $\mult(f)$ such that 
\[
\{I(f,h); \ h\in \mathcal{F}\}= S(f)\setminus \ (S(f)+c(f)-1)
\]
is a representative set of generators of the $k$-vector space $\R/J(f)$. 
\end{theorem}

We postpone the proof of this theorem until the next section, since it is long and quite technical.

In order to use Theorem \ref{theorem1} to prove Theorem \ref{tame}, we will need a kind of Weierstrass Preparation Theorem, which in our case is not suitable, since it makes use of multiplication by units that affects Milnor's number. The solution is given by the Levinson Preparation Theorem which was originally proved in \cite{Le} over $\mathbb C$, but may be adapted without major changes in order to work in arbitrary characteristic. For the proof we refer to the original paper \cite{Le}. \medskip

\noindent {\bf Theorem} \ (Levinson's Preparation Theorem). \label{levinson}
{\em Let $f(X,Y)\in k[[X,Y]]$, where $k$ is an algebraically closed field of arbitrary characteristic $p$. Suppose that $f$ contains for some integer $n>1$ a monomial of the form    $aY^n$ and let $n$ be minimal with this property. If $p \nmid n$, then there exists an automorphism $\psi$ of $k[[X,Y]]$ such that
\[
\psi(f)=A_n(X)Y^n+A_{n-1}(X)Y^{n-1}+\cdots+A_1(X)Y+A_0(X),\] where
$A_i(X)\in k[[X]]$, for all $i$, $A_{n-1}(0)=\cdots=A_1(0)=A_0(0)=0$ and $A_n(0)\neq 0.$
}\medskip

\begin{corollary} Let $f\in k[[X,Y]]$ be irreducible where $k$ is algebraically closed of characteristic $p$. If $p\nmid \mult(f)$, then there exists an automorphism $\varphi$ of $k[[X,Y]]$ such that 
\[\varphi(f)=Y^n+B_{n-1}(X)Y^{n-1}+\cdots+ B_{1}(X)Y+B_0(X),\]
where
$B_i(X)\in k[[X]]$ and $\mult(B_{n-i})>i$, for all $i=1,\ldots,n$.
\end{corollary}
\pf Since $f$ is irreducible, we have that $f=L^n+hot$, where $L$ is a linear form in $X$ and $Y$. By changing coordinates, we may assume that $f$ is as in the Levinson Preparation Theorem.  Now, since $p\nmid n$, we take an $n$-th root of $A_n(X)$ and perform the change of coordinates $Y\mapsto YA_n^{\frac{1}{n}}$ and $X\mapsto X$. So, after only changes of coordinates $\varphi$, we have that
\[\varphi(f)=Y^n+B_{1}(X)Y^{n-1}+\cdots+ B_{n-1}(X)Y+B_n(X),\]
is a Weierstrass polynomial, that is, $\mult(B_{n-i}(X)) > i$, for $i=1,\ldots, n$.
\cqd

\noindent{\bf Proof of Theorem \ref{tame}:} From Deligne's results in \cite{De} one always has   $\mu(f)\geqslant c(f)$.

Now, after a change of coordinates, that do not affect the result, we may assume that $f$ is a Weierstrass polynomial. For every $\alpha \in S(f)\setminus (S(f)+c(f)-1)$, take an element $g\in k[[X,Y]]$ such that $I(f,g)=\alpha$ and after dividing it by $f$ by means of the Weierstrass Division Theorem, we get in this way a family $\mathcal F$ as in Theorem \ref{theorem1}.

Theorem \ref{theorem1} asserts that the residue classes of the elements in $\mathcal{F}$ generate $k[[X,Y]]/J(f)$, hence $\mu(f)\leqslant \# \left(S(f)\setminus (S(f)+c(f)-1)\right)$.
The result will then follow from the next lemma that asserts that the number in the right hand side of the inequality is just $c(f)$. 

The $\mu$-stability follows from the fact that for every invertible element $u$ in $k[[X,Y]]$,
both power series $f$ and $uf$ can be individually prepared to
Weierstrass form by means of a change of coordinates that does not alter the semigroup, nor the Milnor numbers. Hence, $\mu(f)=c(f)=\mu(uf).$ \cqd

\begin{lemma} \label{contagem}
 $\# \left(S(f)\setminus (S(f)+c(f)-1)\right)=c(f)$.
\end{lemma}
\pf In fact, to every $i\in\{0,1,\ldots,c(f)-1\}$ we associate
$s_i\in S(f)\setminus (S(f)+c(f)-1)$ in the following way:
\[
s_i=\left\{\begin{array}{ll}i, & \mbox{if} \ i\in S(f) \\ i+c(f)-1, & \mbox{if}
\ i\not\in S(f). \end{array}\right.\]

The map $i\mapsto s_i$ is injective since $S(f)$ is a symmetric semigroup. On the other hand,   the map is surjective, because, given $j\in S(f)\setminus (S(f)+c(f)-1)$, we have $j=s_j$ if  $j\leqslant c(f)-1$; otherwise, if $j=i+c(f)-1$ for some $i>0$, then again by the symmetry of $S(f)$, it follows that $j$ does not belong to $S(f)$ and therefore $j=s_i$.
\cqd

We believe that the converse of Theorem \ref{tame} is true, in the sense that if $\mu(f)=c(f)$, then $S(f)$ is a tame semigroup, or, equivalently, if $p$ divides any of the minimal generators of $S(f)$, then $\mu(f)>c(f)$. If this is so, we would conclude from our result that if $\mu(f)=c(f)$, then $f$ is $\mu$-stable.

To reinforce our conjecture, observe that the result of \cite{GB-P} proves it when $\mult(f)<p$. The following example is a situation where the converse holds and is not covered by the result in \cite{GB-P}.

\begin{example} Let $p$ be any prime number and $n$ and $m$ two relatively prime natural numbers such that $p\nmid n$, then all curves given by $f(X,Y)=Y^n-X^{mp}$ do not satisfy the condition $\mu(f)=c(f)$, since $\mu(f)=\infty$ and $c(f)=(n-1)(mp-1)$. So, for all $p<n$, we have examples for the converse of our result not covered by \cite{GB-P}.
\end{example}

Anyway, the other possible converse of \ref{tame}, namely, if $f$ is $\mu$-stable then $S(f)$ is tame, is not true, as one may see from the following example.

\begin{example} Let $f=(Y^2-X^3+X^2Y)^2-X^{11}Y\in k[[X,Y]]$, where $\car{k}=5$. Since $f^3\in \M T(f)^3$ (verified with Singular), then $f$ is $\mu$-stable, but its semigroup of values $S(f)=\langle 4,6,25\rangle$ is not tame.
\end{example}

Finally, we observe that the fundamental result used in \cite{GB-P} to prove that $c(f)=\mu(f)$ if and only if $S(f)$ is tame, under the assumption that $p>\mult(f)$, was Lemma 3.3 in that paper that asserts that, in this situation, one has that $I(f,f_Y)\geqslant \mu(f)+\mult(f)-1$, with equality if and only if $S(f)$ is tame. The above inequality is false if one does not assume that $p>\mult(f)$, as we show in the following example.

\begin{example} Let $f=(Y^{9}-X^{13})^2-X^{12}Y\in k[[X,Y]]$. In this case, we have $S(f)=\langle 18,26,301\rangle$, so $c(f)=492$. If  $p=13$, then $p<\mult(f)$ and $\mu(f)=559$  (computed with Singular \cite{DGPS}), hence
\[
I(f,f_Y)=\mult(f)+c(f)-1 =18+492-1= 509 <  576  = \mu(f)+\mult(f)-1.
\]
\end{example}

\section{Proof of Theorem 5.5}

We start with an auxiliary result. Let $f\in \R$ be an irreducible Weierstrass polynomial in $Y$ of degree $n=v_0$, where $S(f)=\langle v_0,\ldots, v_g\rangle$, and $I(f,X)= v_0$ and $I(f,Y)= v_1$. Consider the function defined by $v(h)\colon=I(f,h)$ for $h\in \R\setminus \langle f \rangle$. Consider also the $k[[X]]$-submodule $V_{n-1}$ of $k[[X,Y]]$ generated by $1,Y,\ldots,Y^{n-1}$,
and let $h_0=1, h_1, \ldots, h_{n-1}$ be polynomials in $Y$ such that 
\[ V_{n-1} =  k[[X]]\oplus k[[X]]h_1 \oplus \cdots \oplus k[[X]]h_{n-1},\]
and their residual classes $y_i$ are the Ap\'ery generators of ${\mathcal O}_f$ as a free $k[[X]]$-module (cf. \cite{He} Proposition 6.18).

The natural numbers $a_i=v(y_i)$, $i=0,\ldots,n-1$, form the Ap\'ery sequence of $S(f)$, so they are such that $0=a_0<a_1< \cdots <a_{n-1}$ and $a_i \not\equiv a_j \bmod n$ for $i\neq j$ (cf. \cite{He} Proposition 6.21).

We have the following result.

\begin{proposition}\label{apery}
 Let $I$ be an $\mathfrak{m}$-primary ideal of $\R$ and $h\in V_{n-1}.$
If $v(h)>>0$ then $h\in I$.
\end{proposition}
 
\pf Since the ideal $I$ is $\mathfrak{m}$-primary, there exists a natural number $l$ such that
$\mathfrak{m}^l\subset I$.

Now, write $h=b_0+b_1h_1+\cdots+b_{n-1}h_{n-1}$,
with $b_i\in k[[X]]$, for all $i$. Since $v(b_i) \equiv 0 \bmod n$, $v(h_i)=a_i$ and $a_i\not \equiv a_j \bmod n$, for $i,j=0, \ldots,n-1$, with $i\neq j$, we have that 
\[v(h)=\min_i\{v(b_i)+a_i\}\leqslant\min_i\{v(b_i)\}+a_{n-1}.\]
Hence, $v(h)>>0$ implies that for a given natural number $l$ we have that $\min_j\{v(b_j)\}>l$, hence, $h\in \mathfrak{m}^l\subset I$, as we wanted to show.\cqd

We will need another auxiliary result that appears in \cite{Ca} (Proposition 7.4.1), under the name \textit{Delgado's Formula}, proved over $\C$, which we extend to arbitrary algebraically closed fields.

\begin{lemma}
Let $p=\car{k}$ and $f, g\in k[[X,Y]]$ be a non invertible power series with $f$ irreducible.
Assume that $v_0=\mult(f)$ is not divisible by $p$ and define
$[f,g]=f_Xg_Y-f_Yg_X$. Then one has \[I(f,[f,g])\geqslant I(f,f_Y)-I(f,X)+I(f,g),\]
with equality if and only if $p$ does not divide $I(f,g)$.
\end{lemma}
\pf
Since either $x$ or $y$ is a transversal parameter for $f=0$, we may assume without loss of generality that $x$ is a transversal parameter (the proof in the other case is similar). Let  $(x(t),y(t))$ be a parametrization of $f=0$.  Since $p$ does not divide $v_0$ we have that $\ord_t(x')=v_0-1$ and $\ord_t(x')\leqslant\ord_t(y')$. This shows in particular that $\dfrac{y'}{x'}\in k[[t]]$.\\
\\
Also, $f(x(t),y(t))=0$ implies $f_X(x(t),y(t))x'+f_Y(x(t),y(t))y'=0$, hence
\[
\begin{array}{rcl}
(f_Xg_Y-f_Yg_X)(x(t),y(t)) & = &
 -f_Y(x(t),y(t))\dfrac{y'}{x'}g_Y(x(t),y(t))-f_Y(x(t),y(t))g_X(x(t),y(t))\\
 & = & -\dfrac{1}{x'}f_Y(x(t),y(t))(g'(x(t),y(t))).
\end{array}
\]

On the other hand, since $I(f,g)=\ord_t\Big(g(x(t),y(t))\Big)$, we have
\[I(f,g)-1\leq\ord_t\Big(g'(x(t),y(t))\Big),\] with equality if and only if $p\nmid I(f,g)$.
 It follows that
\begin{eqnarray*}I(f,[f,g])
& = & \ord_t\Big(\dfrac{1}{x'}\Big)+I(f,f_Y)+\ord_t\Big( g'(x(t),y(t)) \Big)\\
& \geqslant & (1-v_0)+I(f,f_Y)+I(f,g)-1\\
& = & I(f,f_Y)-I(f,X)+I(f,g),
\end{eqnarray*}
where equality holds if and only if $p\nmid I(f,g).$\cqd

\begin{remark} \label{delgado}
If $f$ is a Weierstrass polynomial in $Y$ of degree $v_0$ and $p\nmid v_0$, it is well known that
$I(f,f_Y)=c(f)+v_0-1$ (cf. \cite{Za}, proved over $\C$, but same proof works under our hypothesis). So, we conclude that 
$$I(f,[f,g])\geqslant I(f,g)+c(f)-1,$$ 
with equality if and only if $p\nmid I(f,g)$.
\end{remark}

Now, under the assumptions that $f$ is a Weierstrass polynomial in $Y$ with $S(f)=\langle v_0,\ldots, v_g\rangle$ and $p\nmid v_0$, one may associate the Abhyankar-Moh approximate roots  (cf. \cite{A-M}), which are irreducible Weierstrass polynomials $f_{-1}=X,f_0=Y,\ldots,f_g=f$ such that, for each $j$, one has $S(f_j)= \langle \frac{v_0}{e_{j}}, \ldots, \frac{v_j}{e_{j}} \rangle$ and $I(f,f_j)=v_{j+1}$, satisfying a relation, where $\deg$ stands for degree as polynomial in $Y$, 
\[f_j=f^{n_j}_{j-1}-\sum_{i=0}^{n_j-2}a_{ij}f_{j-1}^i,\]
where $a_{ij}$ are polynomials in $Y$ of degree less than $\deg (f_j)=v_0/e_j$ for $j=-1,\ldots,g-1$.

So, from Remark \ref{delgado} we have that
\begin{equation} \label{cordelgado}
I(f,[f,f_{j-1}])\geqslant v_j+c(f)-1, \ \ \hbox{with equality if and only if} \  p\nmid v_j.
\end{equation} 
This implies that if $p\nmid v_0v_1\cdots v_g$, then $S(f)^*+c(f)-1\subset v(J(f))$, where $S(f)^*=S(f)\setminus \{0\}$.

The key result to prove Theorem \ref{theorem1} is Proposition \ref{degreereduction} below that 
will allow us to construct elements in $J(f)\cap V_{n-1}$ whose intersection multiplicity with $f$ sweep the set $S(f)^*+c(f)-1$.

\begin{proposition}\label{degreereduction} Let $f\in k[[X,Y]]$ be a Weierstrass polynomial in $Y$ of degree $v_0$, where $k$ is an algebraically closed field of characteristic 
$p\geqslant 0$. Let $S(f)=\langle v_0,\ldots,v_g\rangle$ and suppose that $p\nmid v_0v_1\cdots  v_g$. Given $s\in S(f)^*$, there exists $q_{s}\in J(f)$, polynomial in $Y$, satisfying
\begin{enumerate}[{\rm (i)}]
\item $\deg {q_{s}}<\deg {f}=v_0$;
\item $I(f,q_{s})=s+c(f)-1$.
\end{enumerate}
\end{proposition}
\pf The proof will be by induction on the genus $g$ of $f$. We will construct step by step the polynomial  $q_s$ which will be of the form $q_s=q_{f,s}=\sum_i P_i[f,f_{j_i}]$ (an infinite sum, possibly) where each $f_{j_i}$ is an approximate root of $f$ and the $P_i$ are monomials in the approximate roots of $f$ satisfying the following conditions:
\begin{equation} \label{condition}
\left\{\begin{array}{ll} I(f,P_1f_{j_1})=s; & \phantom{-}\\ I(f,P_if_{j_i})>s, & \mbox{if}
\ i\neq 1;\\ \deg {f_{j_i}P_i}<\deg {f}, & \mbox{ for all } i.
\end{array} \right.
\end{equation}

If $g=0$, we have $f=Y$, so $J(f)=\R$. Given
$s\in\mathbb{N}=S(f)^*$, set
\[q_{f,s}:=X^{s-1}[f,X].\]
It is easy to check that $q_{f,s}$ satisfies (\ref{condition}) and the conclusion of the proposition.

Inductively, we assume that the construction was carried on for branches of genus $g-1$.
Consider the approximate root $f_{g-1}$ of $f$ of genus $g-1$.
Since $e_{g-1}=n_g$ and $n_gv_g\in \langle v_0,\ldots, v_{g-1} \rangle$, we have
\[S(f)= \langle v_0,\dots,v_g\rangle \subset \Big\langle\frac{v_0}{n_{g}},\dots,
\frac{v_{g-1}}{n_{g}}\Big\rangle =S(f_{g-1}).\] 

For $t\in S(f_{g-1})^*$, the inductive
hypothesis guarantees the existence of a $Y$-polynomial
\[q_{f_{g-1},t}=\sum_i P_i[f_{g-1},f_{j_i}],\] 
where each $f_{j_i}$ is one of the approximate roots $f_{-1},f_0,\ldots,f_{g-2}$ and $P_i$ are monomials in these approximate roots satisfying (\ref{condition}) and the conclusion of the proposition, with $f_{g-1}$ and $v_0/e_{g-1}$ replacing $f$ and $v_0$, respectively. Using this $q_{f_{g-1},t}$ we introduce the following \textit{auxiliary polynomial}
\[\tilde{q}_{f_{g-1},t}:=\sum_i P_i[f,f_{j_i}].\]
To begin with, we will estimate the degree in $Y$ of these polynomials.
The inductive hypothesis gives $\deg {q_{f_{g-1},t}}<\deg {f_{g-1}}$ and
$\deg P_i\leqslant \deg {P_if_{j_i}}<\deg {f_{g-1}}$, for all $i$. On the other hand, the Abhyankar-Moh's relation $f=f_{g-1}^{n_{g}}-G$, where $G=a_{n_{g}-2}f_{g-1}^{n_{g}-2}+\dots+a_0$
and $\deg {a_{i}}<\deg {f_{g-1}}$, gives the inequality $\deg {G} < (n_g-1)\deg {f_{g-1}} = \deg {f}-\deg {f_{g-1}}.$ We also have 
$\deg {[G,f_{j_i}]}=\deg {(G_Xf_{j_i,Y}-G_Yf_{j_i,X}})\leqslant\deg {G}+\deg {f_{j_i}}-1.$
Therefore,
\begin{eqnarray*}\deg{P_i[G,f_{j_i}]}  &\leqslant& \deg{P_i}+\deg{G}+\deg{f_{j_i}}-1\\
 & < & \deg {P_i}f_{j_i}+ \deg {f}-\deg {f_{g-1}} -1\\
 & < & \deg {f_{g-1}}+\deg {f}-\deg {f_{g-1}}=\deg {f},
\end{eqnarray*}
which together with the identity
\[\tilde{q}_{f_{g-1},t}=\sum_i P_i[f_{g-1}^{n_g}-G,f_{j_i}]=
n_gf_{g-1}^{n_g-1}q_{f_{g-1},t}-\sum_i P_i[G,f_{j_i}],\]
give the estimate
$$\deg{\tilde{q}_{f_{g-1},t}}<\deg{f}, \ \ \forall t\in S(f_{g-1})^*.$$\smallskip

\noindent\textbf{Claim $1$:} \textit{For $t\in S(f_{g-1})^*$ we have
$I(f,\tilde{q}_{f_{g-1},t})=c(f)-1+n_gt.$}\\
\\
Indeed, since no generator of $S(f)$ is multiple of $p$, we have from (\ref{cordelgado})
\begin{eqnarray*} 
I(f,P_i[f,f_{j_i}])&=& I(f,P_i)+I(f,[f,f_{j_i}])\\
                   &=& I(f,P_i) + I(f,f_{j_i})+c(f)-1 \\
                   &=& I(f,P_if_{j_i}) +c(f)-1.
                   \end{eqnarray*}
On the other hand, since the $P_if_{j_i}$ are products of approximate roots of $f_{g-1}$ (so, also of $f$), and $I(f_{g-1},P_1f_{j_1})=t$, it follows that $I(f,P_1f_{j_1})=n_gt$. Now, since from (\ref{condition}), the intersection number $I(f,P_if_{j_i})$ assumes its minimum value once for $i=1$, when it is equal to $t$, we have
\begin{eqnarray*}I(f,\tilde{q}_{f_{g-1},t})&=&
I(f,\displaystyle\sum_i P_i[f,f_{j_i}])\\
&=& I(f,P_1f_{j_1}) +c(f)-1\\
& = & n_gt+c(f)-1.
\end{eqnarray*}
\cqd

The family of polynomials
$\{ \tilde{q}_{f_{g-1},t}; \ t \in S(f_{g-1})^* \}$
just introduced will be used in the construction of the family 
$\{ q_{f,s}; \ \ s\in S(f)^*\}$ as announced in the proposition.

To this end, observe that each element $s$ of $S(f)^*$ decomposes
uniquely as 
$$s=n_gt+wv_g, \ \ \hbox{with} \ t\in S(f_{g-1}), \ w\in \{0,1,\dots,n_g-1\}.$$

Now, we break up the analysis in three cases.\\
\\
{\bf Case $1$:} $s=n_gt$. From Claim $1$, we have
$$s+c(f)-1=n_gt+c(f)-1=I(f,\tilde{q}_{f_{g-1},t}).$$
The estimate on the degree of $\tilde{q}_{f_{g-1},t}$, made just before Claim $1$, allows us to deduce that the series
\[q_{f,s}:=\tilde{q}_{f_{g-1},t}\]
has all the required properties, which proves the proposition in this case. \medskip

\noindent {\bf Case $2$:} $s=v_g$. In this case $q_{f,v_g}:=[f,f_{g-1}]$
works because, since $p\nmid v_g$, we have from Remark \ref{delgado} that
$I(f,q_{f,v_g})=v_g+c(f)-1.$ Moreover, using the preceding notations and estimates we get
\begin{eqnarray*}\deg{q_{f,v_g}} & = & \deg{[f_{g-1}^{n_g}-G,f_{g-1}]}\\
 & = & \deg {[f_{g-1},G]}\\
 & \leqslant& \deg{G}+\deg{f_{g-1}}-1\\
 & < & (\deg{f}-\deg{f_{g-1}})+\deg{f_{g-1}}-1\\
 & < & \deg{f}.\\
\end{eqnarray*}

\noindent {\bf Case $3$:} $s>v_g$ and $w>0$. Notice that from the conductor formula one gets that
$c(f)-1=n_g(c(f_{g-1})-1)+(n_g-1)v_g$, and since $n_iv_i<v_{i+1}$, it follows that  $s>v_g>n_g(c(f_{g-1})-1)$. On the other hand, since gcd$(v_g,n_g)=1$, we have that $n_g\nmid s$.

The proposition, in this case, will be established by using the following result that gives a method to reduce degree while preserving intersection multiplicities with $f$ and residual classes modulo $J(f)$. \medskip

\noindent \textbf{Claim $2$:} \textit{Let $s\in\mathbb{N}^*$ be such
that $n_g\nmid s$ and $s>n_g(c(f_{g-1})-1)$. Suppose that we have a
$Y$-polynomial $h$ such that 
$$\deg{h}<\deg{f} \ \ \hbox{and} \ \ I(f,h)=c(f)-1+s,$$ then there  exists a $Y$-polynomial $h'$, such that 
$$\deg{h'} < \deg{f}-\deg{f_{g-1}}, \ \ \ \ I(f,h')=I(f,h),$$
and  
$$h-h'=\displaystyle\sum_{j}\alpha_j\tilde{q}_{f_{g-1},u_j} , \ \  \alpha_j\in k, \ \ u_j\in  S(f_{g-1}), \ \ n_gu_j>s, \forall j.$$}

Indeed we begin by dividing $h$ by $f_{g-1}^{n_g-1}$. Then we get
$h=f_{g-1}^{n_g-1}h_0''+h_0'$ where $\deg{h_0'}<\deg{f_{g-1}^{n_g-1}}=
\deg{f}-\deg{f_{g-1}}.$ The rough idea of the proof is to eliminate the term
$f_{g-1}^{n_g-1}h_0''$ in the preceding relation using the polynomials
$\tilde{q}_{f_{g-1},u}$ where $u\in S(f_{g-1})^*$. This will be done iteratively, in possibly infinitely many steps, with the help of the following auxiliary result.\\

\noindent \textbf{Claim $3$:}  \textit{
With the same conditions as above, we have
$I(f,h_0'')=n_gI(f_{g-1},h_0'')$  and
$I(f,h_0''f_{g-1}^{n_g-1})\neq I(f,h_0').$}

We will prove this claim after the conclusion of the proof of Claim $2$, given below.

Using the formula $c(f)-1=n_g(c(f_{g-1})-1)+(n_g-1)v_g$ and Claim $3$, we get
\[I(f,h_0''f_{g-1}^{n_g-1})-(c(f)-1)=n_g[I(f_{g-1},h_0'')-c(f_{g-1})+1].\]
On the other hand, since $I(f,h)-(c(f)-1)=s$ and $n_g\nmid s$, it follows that 
 \[I(f,h_0')=I(f,h)<I(f,f_{g-1}^{n_g-1}h_0'').\]

So, from the first part of Claim $3$ and the above inequality, we get  
$$n_gI(f_{g-1},h_0'')=I(f,h_0'')>I(f,h)-I(f,f_{g-1}^{n_g-1})=
s+c(f)-1-(n_g-1)v_g.$$

Defining $u_1=I(f_{g-1},h_0'')-c(f_{g-1})+1$, it follows that
$$c(f_{g-1})-1+u_1= I(f_{g-1},h_0'')>\frac{s}{n_g}+c(f_{g-1})-1>2(c(f_{g-1})-1),$$
allowing us to conclude that $u_1\in S(f_{g-1})^*$. 

The inductive hypothesis guarantees the existence
of a polynomial $q_{f_{g-1},u_1}$ satisfying all requirements in (\ref{condition}) and the conclusion in Proposition \ref{degreereduction}.

From Claim $1$, we have
\[I(f,\tilde{q}_{f_{g-1},u_1})=c(f)-1+n_gu_1=I(f,h_0''f_{g-1}^{n_g-1}).\]
So, after multiplication by a suitable $\alpha_1\in k^{*}$, we get that $h_1=h_0''f_{g-1}^{n_g-1}-\alpha_1\tilde{q}_{f_{g-1},u_1}$ satisfies the inequality
\begin{equation} \label{firststep} I(f,h_1)> I(f,h_0''f_{g-1}^{n_g-1})>I(f,h) \ (=c(f)-1+s).\end{equation}

This allows us to write 
\[
h=h_1+\alpha_1 \tilde{q}_{f_{g-1},u_1}+h_0', \ \ \hbox{with} \ \ I(f,h_1)>I(f,h) \ \  \hbox{and} \ \ I(f,h_0')=I(f,h).
\]

From (\ref{firststep}) we have that there exists $s_1\in {\mathbb N}^*$ such that
\[
I(f,h_1)=c(f)-1+s_1 >c(f)-1 +n_gu_1>c(f)-1+s.
\]
So, $s_1>s$ and $n_gu_1>s$.

In the next step we proceed differently according to the divisibility of $s_1$ by $n_g$.
Suppose $n_g\mid s_1$, say $s_1=n_gu_2$. In this case, by the above inequality we have
\[2(c(f_{g-1})-1)<c(f_{g-1})-1+u_1<c(f_{g-1})-1+u_2.\] So, it follows that
$u_2\in S(f_{g-1})^*$. Hence, there exists a polynomial
$q_{f_{g-1},u_2}$ such that
\[I(f,h_1)=I(f,\tilde{q}_{f_{g-1},u_2} ) \] and again we may choose $\alpha_2\in k$ in such a way that if $h_2=h_1-\alpha_2\tilde{q}_{f_{g-1},u_2}$, we have
$I(f,h_2)>I(f,h_1)$.
Hence, we get $h=h_2+\alpha_1\tilde{q}_{f_{g-1},u_1}+h_0'+\alpha_2\tilde{q}_{f_{g-1},u_2}+h_1'$, where $h_1'=0$, in this case. Notice that $n_gu_2>n_gu_1>s$.

If, however, $n_g\nmid s_1$, we are in position to repeat the preceding
process of division by $f_{g-1}^{n_g-1}$ using, this time, $h_1$ instead of $h$.
So $h_1=f_{g-1}^{n_g-1}h_1''+h_1'$. Again,  we deduce that there exist $\alpha_2\in k$ and $u_2\in S(f_{g-1})^*$, with $n_gu_2>s_1>s$, such that if we define $h_2=h_1''f_{g-1}^{n_g-1}-\alpha_2\tilde{q}_{f_{g-1},u_2}$, then we have
$$I(f,h_2)>I(f,h_1)>I(f,h).$$

So, by repeating this process we get  
$$h=h_j+ \sum_{i=1}^j \alpha_i\tilde{q}_{f_{g-1},u_i}+ \sum_{i=0}^{j-1}h_i',$$
with $I(f,h_i')<I(f,h_{i+1}')$, if $h_i'\neq 0$, and $I(f,h_i)<I(f,h_{i+1})$.
Since all power series appearing in the above sum have degree less than $\deg{f}$, it follows, in view of Proposition \ref{apery}, that $h_j\rightarrow 0$ in the
$\mathfrak{m}$-adic topology of $\R$ and the family $\{h_i'\}_{i\in {\mathbb N}}$ is summable. Taking $h'=\sum_j h_j'$ we get Claim $2$.
\cqd

Finally it remains to prove Claim $3$. If $f$ is any irreducible Weierstrass polynomial of degree $n$, then it is easy to see from Proposition \ref{apery} that the set $V_{n-1}$ of all polynomials in $Y$ of degree less than $n$ with coefficients in $k[[X]]$ is a free $k[[X]]$-module with basis 
$$\left\{f^J=f_0^{j_0}f_1^{j_1}\cdots f_{g-1}^{j_{g-1}}; \ J=(j_0,\ldots,j_{g-1}), \ 0\leqslant j_i<n_{i+1}, \ i=0,\ldots,g-1 \right\}.$$ So, every element $h\in V_{n-1}$ may be written uniquely as
\[h=\sum_J a_J(X)f^J = 
f_{g-1}^{n_g-1}h''+h', \ \ \ \ a_J(X)\in k[[X]],
\]
where 
\[
h''=\displaystyle\sum_{j_{g-1}=n_g-1}a_J(X)f_0^{j_0}
\cdots f_{g-2}^{j_{g-2}}, \ \ \ \ h'=\displaystyle\sum_{j_{g-1}\leqslant
n_{g}-2}a_J(X)f_0^{j_0}\cdots f_{g-1}^{j_{g-1}}.
\]

First of all we will check that $I(f,h')\neq I(f,f_{g-1}^{n_g-1}h'')$. In fact, in $h'$ there is a unique term such that
\[I(f,h')=I(f,a_J(X)f_0^{j_0}\cdots f_{g-1}^{j_{g-1}})=\sum_{i=-1}^{g-1}j_iv_{i+1},\]
where $j_{-1}=\ord_Xa_J(X)$. Also, in $f_{g-1}^{n_g-1}h''$ there is a unique term satisfying
\[I(f,f_{g-1}^{n_g-1}h'')=I(f,a_L(X)f_0^{l_0}\cdots f_{g-2}^{l_{g-2}}f_{g-1}^{n_g-1})=
\sum_{i=-1}^{g-2}l_iv_{i+1}+(n_g-1)v_g, \]
where $ l_{-1}=\ord_X(a_L(X))$.

Since each element in $S(f)$ is written in a unique way as $\sum_{i=-1}^{g-1}j_iv_i$ with 
$j_{-1}\in {\mathbb N}$ and $0\leqslant j_i\leqslant n_{i+1}-1$, the inequality follows.

Also, it is clear from the way we wrote $h''$ that $I(f,h'')=n_{g}I(f_{g-1},h'')\in n_gS(f_{g-1})$.

Now, to finish the proof of Claim $3$ we only need to check that $h''$
and $h'$ are indeed the quotient and the remainder, respectively, of the division of
$h$ by $f_{g-1}^{n_g-1}$. We will do this by estimating the degree of $h'$ and , hence,
conclude by the uniqueness of the remainder and the quotient in the euclidean algorithm.
Indeed, for every summand in $h'$ we have
$\deg{(a_J(X)f_0^{j_0}\cdots f_{g-1}^{j_{g-1}})}< \deg{f_{g-1}^{n_g-1}}$, which shows that  $\deg{h'}<\deg{f_{g-1}^{n_g-1}}$.\cqd \smallskip

Now we return to the construction of the polynomial $q_{f,s}$ in the remaining
Case $3$, that is, when $s=n_gt+v_gw$ with $s>v_g$ and $w>0$. 

Observe that if $t=0$ and $w=1$, then from Case $2$ we have $q_{f,v_g}=[f,f_{g-1}]$. Now, we apply Claim $2$ to $h=q_{f,v_g}$ in order to find
$h'=(q_{f,v_g})'$ with degree less than $\deg{f}-\deg{f_{g-1}}$ satisfying
\[I(f,(q_{f,v_g})')=I(f,q_{v_g,f})=c(f)-1+v_g\] and
\[ (q_{f,v_g})'=q_{f,v_g}+\displaystyle\sum_{j}\alpha_j
\tilde{q}_{f_{g-1},u_j}.\]
Using this, we define $ q_{f,2v_g}:=f_{g-1}(q_{f,v_g})'$. Clearly, we have 
$\deg{q_{f,2v_g}}<\deg{f}$ and $I(f,q_{f,2v_g})=c(f)-1+2v_g$.
Hence, it remains to show that $q_{f,2v_g}$ satisfies (\ref{condition}) in order to make possible our inductive process.

We have \[q_{f,2v_g}=f_{g-1}[f,f_{g-1}]+
\displaystyle\sum_{j}\alpha_jf_{g-1}\tilde{q}_{f_{g-1},u_j}=\sum_iP_i[f,f_{j_i}].\]
This shows that $q_{f,2v_g}$ has the required format.
Finally, we need to check the statement about intersection indices. We are going to show
that $P_1=f_{j_1}=f_{g-1}$. In order to do so, it is enough to show that
for each index $j$ in the above sum, the polynomial
\[\tilde{q}_{f_{g-1},u_j}=\sum_l P_l'[f,f_{j_l}],\] where $f_{j_l}$ is one of the approximate roots $f_{-1},f_0,\ldots,f_{g-2}$ and the $P_l'$ are monomials in these approximate roots, 
is such that
$I(f,f_{g-1}P_i'f_{j_i})>2v_g$. Indeed, from the inductive hypothesis we have
\[
I(f,f_{g-1}P_l'f_{j_l})  =  v_g+I(f,P_l'f_{j_l})
 =  v_g+n_gI(f_{g-1},P_l'f_{j_l})> v_g+n_gu_1> 2v_g,
 \]
where the last strict inequality is justified by the fact that from Claim $2$ one has $n_gu_1>v_g$.

We apply again Claim $2$ to obtain $(q_{f,2v_g})'$ which multiplied by $f_{g-1}$ produces $q_{f,3v_g}$. Now, we repeat this procedure until we get the
polynomial $q_{f,wv_g}=f_{g-1}(q_{f,(w-1)v_g})'$, satisfying the
proposition for $wv_g\in S(f)^*$. Since $s=n_gt+wv_g$, we consider the
polynomial $q_{f_{g-1},t}=\sum P_i''[f_{g-1},f_{m_i}]$ and collect
$P_1'', \ f_{m_1}$ so that $I(f_{g-1}, P_1''f_{m_1})=t.$
Finally, define \[q_{f,s}:=P_1''f_{m_1}(q_{f,wv_g})'.\]
It is now immediate to verify that $q_{f,s}$ satisfies (\ref{condition}) and the conclusion of the proposition, finishing its proof. \cqd

With these tools at hands, we may conclude the proof of Theorem \ref{theorem1}.\\

\noindent \textit{\bf Proof of the Theorem \ref{theorem1}:} Choose ${\mathcal{F}}$ with minimal number of elements, so from Lemma \ref{contagem} it follows that $\#{\mathcal{F}}=c(f)$. We will show that
the set $\overline{\mathcal{F}}$ generates $\R/J(f)$ as a $k$-vector space.
In particular, this will show also that $\mu(f)\leqslant c(f)$ when $S(f)$ is tame.
In order to do this it is enough to show that there exists a decomposition
$\R=\langle \mathcal{F}\rangle + J(f),$ where $\langle \mathcal{F}\rangle$ denotes the $k$-vector space spanned by the elements of $\mathcal{F}=\{\varphi_1,\ldots,\varphi_{c(f)}\}$. 

Given any element $h\in \R$ we can divide it  by the partial derivative $f_Y$
which, under our assumptions, is a $Y$-polynomial of degree $v_0-1$.
The remainder of the division is a $Y$-polynomial $h'$ of degree less than $v_0-1$
and it is sufficient to show that $h'$ belongs to
$\langle \mathcal{F}\rangle + J(f).$

If $I(f,h')\in S(f)\setminus (S(f)+c(f)-1)$ then, according to the definition of $\mathcal{F}$, there is an element $\varphi_{i_{s_0}}$ such that
$s_0:=I(f,h')=I(f,\varphi_{i_{s_0}})$. Hence, there is a constant $\alpha_{s_0}\in k$ such that
\[I(f,h'-\alpha_{s_0}\varphi_{i_{s_0}})=:s_1>s_0.\] 

If, otherwise $s_0=I(f,h')\in S(f)^*+c(f)-1$, then choose an element $q_{f,s_0}$ in $J(f)$ polynomial in $Y$ of degree less then $\deg{f}$, such that $s_0=I(f,q_{f,s_0})$. Hence, there is a constant $\beta_{s_0}\in k$ such that
\[I(f,h'-\beta_{s_0}q_{f,s_0})=:s_1>s_0.\] 

We carry on this process that increases intersection indices to eventually achieve 
\[s_r=I(f,h'- \sum_s\beta_{s}q_{f,s}+  \sum_s \alpha_s\varphi_{i_s})\in S(f)^*+c(f)-1, \ \ \ \forall r\geqslant N.\] Since the elements in $S(f)^*+c(f)-1$ may be realized as intersections indices of $f$ with elements in $J(f)\cap V_{n-1}$ (cf. Proposition \ref{degreereduction}), we produce an element 
\[h'-\sum_s \alpha_s\varphi_{i_s}-\sum_j \beta_s q_{f,s}\]
whose intersection multiplicity with $f$ is big enough and whose degree is less than $\deg{f}$, hence from
Proposition \ref{apery} it belongs to the Jacobian ideal $J(f)$. \cqd
\bigskip

\noindent Authors addresses: {hefez@mat.uff.br; \ joaohelder@hotmail.com; \ rsalomao@id.uff.br}


                             

\end{document}